\documentstyle[amsmath,amssymb,12pt]{article}

\renewcommand{\author}[1]{\medskip{\Large #1}\par\medskip}
\newcommand{\D}{{\cal D}}
\newcommand{\N}{{\Bbb N}}

\newcommand{\R}{{\Bbb R}}
\newcommand{\la}{\langle}
\newcommand{\ra}{\rangle}

\makeatletter\@addtoreset{equation}{section}\makeatother

\allowdisplaybreaks

\begin{document}

\setcounter{page}{1} \setcounter{section}{0} \thispagestyle{empty}

\newcommand{\rom}[1]{{\rm #1}}

\newtheorem{definition}{Definition}[section]
\newtheorem{remark}{Remark}[section]
\newtheorem{proposition}{Proposition}[section]
\newtheorem{theorem}{Theorem}[section]
\newtheorem{corollary}{Corollary}[section]
\newtheorem{lemma}{Lemma}[section]

\thispagestyle{empty}

\begin{center}{\Large \bf A model of the term structure of interest rates based on
 L\'evy  fields}\end{center}

\author {Sergio Albeverio}
\noindent{\sl Institut f\"{u}r Angewandte Mathematik,
Universit\"{a}t Bonn, Wegelerstr.~6, D-53115 Bonn,\\ Germany; SFB
611, Univ.~Bonn, Germany; CERFIM (Locarno); Acc.\ Arch.\ (USI),
Switzerland;  BiBoS, Univ.\ Bielefeld, Germany}\\[2mm]
 E-mail: \texttt{albeverio@uni-bonn.de}

\author{Eugene Lytvynov}

\noindent{\sl Institut f\"{u}r Angewandte Mathematik,
Universit\"{a}t Bonn, Wegelerstr.~6, D-53115 Bonn, Germany;
 SFB 611, Univ.~Bonn, Germany;
BiBoS, Univ.\ Bielefeld, Germany}\\[2mm] E-mail:
\texttt{lytvynov@wiener.iam.uni-bonn.de}

\author{Andrea Mahnig}

\noindent{\sl Institut f\"{u}r Angewandte Mathematik,
Universit\"{a}t Bonn, Wegelerstr.~6, D-53115 Bonn}\\[2mm]
  E-mail:            \texttt{andrea.mahnig@gmx.de}\vspace{10mm}



\begin{abstract}

\noindent An extension of the Heath--Jarrow--Morton model for the
development of instantaneous forward interest rates with
deterministic coefficients and Gaussian as well as L\'evy field
noise terms is given. In the special case where the L\'evy field
is absent, one recovers a model discussed by D.P.~Kennedy.

\end{abstract}

\vspace{2mm}

 \vspace{2mm}


\section{Introduction}

 Heath, Jarrow, and Morton (1992) (see also Heath et al.\ (1990))
 proposed a model of interest rates and their
associated bond prices in which the price at time $s$ of a bond
paying one unit at time $t\ge s$ is given by \begin{equation}
P_{s,t}=\exp\bigg[-\int_s^t F_{s,u}\,du
\bigg],\label{zwgc}\end{equation} where $F_{s,t,}$, $0\le s\le t$,
is  called the {\it instantaneous forward  rate}, or just the {\it
forward rate}. Then \begin{equation} R_s{:=}F_{s,s},\qquad s\ge
0,\label{one}\end{equation} is called the {\it instantaneous spot
rate}, or just the {\it spot rate}. One also defines the {\it
discounted bond-price process\/} as
\begin{equation}\label{two}Z_{s,t}{:=}P_{s,t}\,\exp\bigg[-\int_0^s R_u\, du\bigg].
\end{equation}

In the Heath--Jarrow--Morton (HJM) model, the forward rates are
supposed to  satisfy the stochastic differential equations
\begin{equation}dF_{s,t}=\alpha(s,t)\,
ds+\sum_{i=1}^m\beta_i(s,t)\, dW_s^i,\label{ucu}\end{equation}
where  $W^1,\dots,W^m$ are independent standard Brownian motions
and $\alpha(s,t)$ and $\beta_i(s,t)$ are processes adapted to the
natural filtration of the Brownian motions. This model was, in
fact, an extension of the earlier work  by Ho et al.\ (1986).

Kennedy (1994) (see also Kennedy (1997)), while following the
approach of modeling the instantaneous forward rates, considered
the case where $\{F_{s,t},\, 0\le s\le t<\infty\}$ is a continuous
Gaussian random field which has independent increments in the
$s$-direction, that is, in the direction of evolution of `real'
time. This framework includes the HJM model in the case where the
coefficients $\alpha(s,t)$ and $\beta_i(s,t)$ in \eqref{ucu} are
deterministic. An important example of application of the Kennedy
model is the case where the forward rates are given by
$F_{s,t}=\mu_{s,t}+X_{s,t}$ with $\mu_{s,t}$ being deterministic
and  $X_{s,t}$  a Brownian sheet (see, e.g., Adler (1981) for this
concept). In the latter case, $F_{s,t}$ has independent increments
also in the $t$ direction. Furthermore, this  may be intuitively
thought of as the situation of \eqref{ucu} driven by an
uncountably infinite number of Brownian motion.  Kennedy (1994)
gave a simple characterization of the discounted bond-price
process to be a martingale. In particular, he showed that the
latter is true if and only if the expectation $\mu_{s,t}$ of
$F_{s,t}$, $0\le s\le t<\infty$, satisfies a simple relation.

In Bj\"ork {et al.}\ (1997a, b)  (see also Bj\"ork  et al.\
(1999)), it was pointed out that, in many cases observed
empirically, the interest rate trajectories do not look like
diffusion processes, but rather as diffusions and jumps, or even
like pure jump processes. Therefore, one needs to introduce a jump
part in the description of interest rates. The authors of these
papers considered the case where the forward rate process
$\{F_{s,t},\,0\le s\le t\}$ is driven by a general marked point
process as well as by a Wiener process (Bj\"ork et al.\ (1997b)),
or by a rather general L\'evy process (Bj\"ork et al.\ (1997a)),
and the maturity time $t\ge 0$ is a continuous parameter of the
model. In particular, an equivalence condition was given for a
given probability measure to be a local martingale measure, i.e.,
for the discounted bond-price process $\{Z_{s,t},\, 0\le s\le t\}$
to be a local martingale for each $t\ge 0$ (Propositions~5.3, 5.5
in Bj\"ork et al.\ (1997a), see also Theorem~3.13 in Heath et al.\
(1992)). This condition, formulated in terms of the coefficients
for the forward rate dynamics, generalizes the result of Heath et
al.\ (1990) which was obtained for the diffusion case.

Other generalizations of the HJM model in which the forward rate
process satisfies  stochastic differential equations with an
infinite number of independent standard Brownian motions (i.e.,
$m=\infty$ in \eqref{ucu}) were proposed in Yalovenko (1998),
Kusuoka (2000), Popovici (2001), and L\"utkebohmert (2002).

In the present paper, we follow the approach of Kennedy (1994,
1997), but suppose that the  forward rates $\{F_{s,t},\,0\le s\le
t<\infty\}$ are given by  a L\'evy field without a diffusion part.
In particular, $\{F_{s,t}\}$ has independent increments in both
$s$ and $t$ directions. Analogously to Kennedy (1994), we give, in
this case, a characterization of the martingale measure. We also
show that, under a slight additional condition on the L\'evy
measure of the field, it is possible to choose the initial term
structure $\{\mu_{0,t},\,t\ge0\}$ in such a way that the forward
interest rates are a.s.\ non-negative. This, of course, was
impossible to reach in the framework of the Gaussian model, which
caused problems  in some situations
 (see Sect.~1 of Kennedy (1997)). We then present two examples of application of our
 results: the cases where $F_{s,t}$ is a ``Poisson sheet'' (this case was discussed in Mahnig (2002)),
 respectively a
 ``gamma sheet.''
 Finally, we mention the
possibility of unification of the approaches of Kennedy and of the
present paper, by considering $F_{s,t}$ as a sum of a Gaussian
field and an independent L\'evy field, and thus having a process
with a  diffusion part as well as  a jump part.

\section{The model based on   L\'evy fields}

Let $\D{:=}C_0^\infty(\R^2)$ denote the space of all real-valued
infinitely differentiable functions on $\R^2$ with compact
support. We equip $\D$ with the standard nuclear space topology,
see, e.g., Berezansky et al.\ (1996). Then $\D$ is densely and
continuously embedded into the real space $L^2(\R^2,dx\, dy)$. Let
$\D'$ denote the dual space of $\D$ with respect to the
``reference''  space $L^2(\R^2,dx\,dy)$, i.e., the dual pairing
between elements of $\D'$ and $\D$ is generated by the scalar
product in $L^2(\R^2,dx\,dy)$. Thus, we get the standard
(Gel'fand) triple $$ \D'\supset L^2(\R^2,dx\,dy)\supset \D.$$ We
denote by $\la\cdot,\cdot\ra$ the dual pairing between elements of
$\D'$ and $\D$. Let  ${\cal C}(\D')$ denote the cylinder
$\sigma$-algebra on $\D'$.

We define a centered L\'evy noise measure as a probability measure
$\nu$ on $(\D',{\cal C}(\D'))$ whose Fourier transform is given by
\begin{equation}\label{1}
\int_{\D'}e^{i\la\omega,\varphi\ra}\,\nu(d\omega)=\exp\bigg[
\int_{\R^2}\int_{\R_+}(e^{i\tau \varphi(x,y)}-1-i\tau
\varphi(x,y))\,\sigma(d\tau )\,dx\,dy\bigg],\qquad \varphi\in\D
\end{equation}
(see, e.g., Ch.~III, Sec.~4 in Gel'fand and Vilenkin (1964)).
Here, $\R_+=(0,+\infty)$ and $\sigma$ is a positive measure on
$(\R_+,{\cal B}(\R_+))$, which is usually called the L\'evy
measure of the process. We suppose that $\sigma$ satisfies the
following condition:
\begin{equation}\label{2} \int_{\R_+}\tau
^2\,\sigma(d\tau)<\infty.
 \end{equation}
The existence of the measure $\nu$ follows from the Minlos
theorem.

For any $\varphi\in\D$, we easily have
\begin{equation}\label{3}\int_{\D'}\la\omega,\varphi\ra^2\,
\nu(d\omega)=\int_{\R_+}\tau ^2\,\sigma(d\tau)
\int_{\R^2}\varphi(x,y)^2\, dx\,dy.\end{equation} Thus, the
mapping $I: L^2(\R^2,dx\,dy)\to L^2(\D',d\nu)$,
$\operatorname{Dom}(I)=\D$, defined by
$$(I\varphi)(\omega){:=}\la\omega,\varphi\ra,\qquad \varphi\in\D,\
\omega\in\D',$$ may be extended by continuity to the whole
$L^2(\R^2, dx\, dy)$. For each $f\in L^2(\R^2,dx\,dy)$,  we set
$\la\cdot,f\ra{:=}If$. Thus, the random variable (r.v.)
$\la\omega, f\ra$ is well-defined for $\nu$-a.e.\ $\omega\in\D'$
and equality \eqref{3} holds for $f$ replacing $\varphi$.

Let $\varkappa:[0,\infty)^2\to[0,\infty)$ be a measurable function
which is locally bounded. For each $s,t\ge0$, we  define the r.v.\
$X_{s,t}$ as follows:
\begin{equation}\label{zdf} X_{s,t}(\omega){:=}\la\omega(x,y),\pmb
1_{[0,s]} (x)\pmb 1_{[0,t]}(y)\varkappa(x,y)\ra,\qquad
\text{$\nu$-a.e.\ $\omega\in\D'$ },\end{equation} where $x,y$
denote the variables in which the dualization is carried out. It
follows from \eqref1 that  $X_{s,t}$ is centered and has
independent increments in both the $s$ and $t$ directions.

We note that, in the  case where $\varkappa(x,y)\equiv 1$,
$\{X_{s,t},\ 0\le s\le t\}$ is a L\'evy process for each fixed
$t>0$. Indeed, it follows from \eqref1 that the Fourier transform
of $X_{s,t}$ is given by $$ \int_{ \D'}e^{i\lambda
X_{s,t}(\omega)}\,\nu(d\omega)=\exp\bigg[st\int_{\R_+} (e^{i\tau
\lambda}-1-i\tau \lambda)\,\sigma(d\tau )\bigg],\qquad
\lambda\in\R.$$ In particular, the L\'evy measure of the process
$\{X_{s,t},\ 0\le s\le t\}$ is equal to $t\sigma$.

Let $F_{s,t}$ be the  forward  rate for date $t$ at time $s$,
$0\le s\le t$. We suppose that
\begin{equation}\label{fwzc} F_{s,t}=\mu_{s,t}+X_{s,t},\qquad 0\le s\le
t,\end{equation} where $\mu_{s,t}$ is deterministic and continuous
in $(s,t)$ on $\{(s,t)\in\R^2: s\le t\}$. The price at time $s$ of
a bond paying one unit at time $t\ge s$ is then given by
\eqref{zwgc}. We note that the random variable $\int_s^t X_{s,u}\,
du$ is $\nu$-a.s.\ well-defined and
\begin{equation}\label{4}\int_s^t X_{s,u}(\omega)\, du=\la
\omega(x,y),\pmb 1_{[0,s]}(x)\pmb
1_{[0,t]}(y)\varkappa(x,y)(t-(s\vee y))\ra\quad\text{for
$\nu$-a.e.\ $\omega\in \D' $ }.\end{equation}  Indeed, for each
$f\in L^2(\R^2,dx\,dy)$, we have by \eqref{3}:
\begin{align}& \int_{\D'} \left(\int_s^t X_{s,u}(\omega)\, du\right)\la
\omega,f\ra\,\nu(d\omega)\notag\\ &\qquad= \int_s^t
\int_{\D'}X_{s,u}(\omega) \la \omega,f\ra\,\nu(d\omega)\,
du\notag\\ &\qquad=\int_{\R_+} \tau ^2\,\sigma(d\tau)\cdot
\int_s^t \int_{\R^2}\pmb 1_{[0,s]}(x)\pmb
1_{[0,u]}(y)\varkappa(x,y)f(x,y)\, dx\, dy\, du\notag\\ &\qquad
=\int_{\R_+} \tau ^2\,\sigma(d\tau)\cdot\int_{\R^2}\pmb
1_{[0,s]}(x)\left(\int_s^t\pmb 1_{[0,u]}(y)\, du
\right)\varkappa(x,y)f(x,y)\, dx\, dy\notag\\ &\qquad =\int_{\R_+}
\tau ^2\,\sigma(d\tau)\cdot \int_{\R^2}\pmb 1_{[0,s]}(x)\pmb
1_{[0,t]}(y)(t-(s\vee y))\varkappa(x,y)f(x,y)\, dx\,
dy,\label{5}\end{align} which implies \eqref{4}.

Analogously to \eqref4, \eqref5, we have: $$ \int_0^s R_u \,du =
\int_0^s \mu_{u,u}\,du+\la \omega(x,y),\pmb 1_{[0,s]}(x)\pmb 1
_{[0,s]}(y)\varkappa(x,y)(s-(x\vee y)), $$ where the spot rate
$R_u$ is defined by \eqref{one}.  Thus,  the discounted bond-price
process $Z_{s,t}$  is well-defined by  \eqref{two}.

 We
 denote by ${\cal F}_s$, $s\ge 0$, the $\sigma$-algebra
generated by the r.v.'s $F_{u,v}$, $0\le u\le s$, $u\le v$, which
describes the information available at time $s$.

\begin{theorem}\label{th1}
The following statements are equivalent\rom:

\noindent \rom{(a)} For each $t\ge 0$\rom, the discounted
bond-price process $\{Z_{s,t},\ {\cal F}_s,\ 0\le s\le t\}$ is a
martingale\rom.

\noindent\rom{(b)} $$ \mu_{s,t}=\mu_{0,t} + \int_0^t \int_0^s
\int_{\R_+}\tau\varkappa(x,y)(1-e^{-\tau\varkappa(x,y)(t-(x\vee
y))})\,\sigma(d\tau)\, dx\, dy $$ for all $s,t\ge0$\rom, $s\le
t$\rom.

\noindent{\rom{(c)}} $$ P_{s,t}={\Bbb E}(e^{-\int_s^t R_u\,
du}\mid {\cal F}_s)\qquad \text{\rom{for all $s,t\ge0$, $s\le t$.
}}$$

\end{theorem}

\noindent{\it Proof}. We first show the equivalence of (a) and
(b). Absolutely analogously to the proof of Theorem~1.1 in Kennedy
(1994), we conclude that (a) is equivalent to the following
condition to hold: \begin{equation}\label{6}
\int_{\D'}\exp\bigg(-\int_{s_1}^t (F_{s_1,u}-F_{s_2,u})\, du
-\int_{s_2}^{s_1} (F_{u,u}-F_{s_2,u})\,du \bigg)\,
d\nu=1\end{equation} for all $0\le s_2\le s_1\le t$. By
\eqref{zdf} and \eqref{fwzc}, \eqref{6} is equivalent to
\begin{align}&\int_{\D'}\exp\bigg( -\int_{s_1} ^t \la \omega(x,y),
\pmb 1_{[s_2,s_1]}(x)\pmb 1_{[0,u]}(y)\varkappa(x,y)\ra\,
du\notag\\ & \qquad\quad\text{}-\int_{s_2}^{s_1} \la \omega(x,y) ,
\pmb 1_{[s_2,u]}(x)\pmb 1_{[0,u]}(y)\varkappa(x,y)\ra\,du
\bigg)\,\nu(d\omega)\notag\\ &\qquad = \exp\bigg( \int_{s_1}^t
(\mu_{s_1,u}-\mu_{s_2,u})\,du+\int_{s_2}^{s_1}(\mu_{u,u}-\mu_{s_2,u})\,du\bigg)\label{ucdg}
\end{align}
 for all $0\le s_2\le s_1\le t$.
Analogously to \eqref4 and  \eqref5, we have
\begin{align*}&\int_{s_1} ^t \la \omega(x,y),
\pmb 1_{[s_2,s_1]}(x)\pmb 1_{[0,u]}(y)\varkappa(x,y)\ra\, du
+\int_{s_2}^{s_1} \la \omega(x,y) , \pmb 1_{[s_2,u]}(x)\pmb
1_{[0,u]}(y)\varkappa(x,y)\ra\,du\\ &\qquad= \la\omega(x,y), \pmb
1_{[s_2,s_1]}(x)\pmb 1_{[0,t]}(y)\varkappa(x,y)(t-(x\vee y))\ra
\qquad \text{for $\nu$-a.e.\ $\omega\in\D'$. }
\end{align*} Hence, it follows from  \eqref1  that
condition  \eqref{ucdg} is equivalent to
\begin{align}&  \int_{\R_+} \int_0^t \int _{s_2}^{s_1}
(e^{-\tau\varkappa(x,y)(t-(x\vee y))}-1+
\tau\varkappa(x,y)(t-(x\vee y))) \, dx\, dy\,
\sigma(d\tau)\notag\\ &\qquad = \int_{s_1}^t
(\mu_{s_1,u}-\mu_{s_2,u})\,du+\int_{s_2}^{s_1}(\mu_{u,u}-\mu_{s_2,u})\,du\label{ige}
\end{align}
 for all $0\le s_2\le s_1\le t$. We remark  that the function under
 the sign of integral on the left hand side of \eqref{ige} is
 integrable. Indeed, let us set $C_t{:=}\sup_{x,y\in[0,t]}\varkappa(x,y)$.
 Then, for all $\tau\in(0,1]$ and $x,y\in[0,t]$, we have
 \begin{align}&|e^{-\tau\varkappa(x,y)(t-(x\vee y))}-1+\tau\varkappa(x,y)
 (t-(x\vee y))|\notag\\
 &\qquad\le \sum_{n=2}^\infty\frac{(\tau C_t t)^n}{n!}\le \tau^2
C_t^2 t^2\exp(\tau C_t t)\le \tau^2 C_t^2 t^2\exp(C_t
t).\label{ibib}\end{align} This, together with the fact that
$\int_{(0,1]}\tau^2\,\sigma(d\tau)<\infty$, yields the
integrability of the integrand on the left hand side of
\eqref{ige} on $(0,1]\times [0,t]\times [s_1,s_2]$. Furthermore,
for $\tau\in(1,+\infty)$ and $x,y\in[0,t]$, we have
$$|e^{-\tau\varkappa(x,y)(t-(x\vee
y))}-1+\tau\varkappa(x,y)(t-(x\vee y))| \le 1+\tau tC_t.$$ This,
together with  the fact that $
\int_{(1,+\infty)}\tau\,\sigma(d\tau)<\infty$, completes the proof
of the integrability of the integrand on  the left hand side of
\eqref{ige}.

We now fix $t>0$ and suppose, for a moment, that $\mu_{s,t}$ has
the following form: \begin{equation}\label{ziguczug}
\mu_{s,t}=\mu_{0,t}+\int_0^t\int_0^s\Psi_t(x,y)\,
dx\,dy,\end{equation} where $\Psi_t(x,y)$ is an integrable
function on $[0,t]^2$. Then,
\begin{equation} \int_{s_1}^t
(\mu_{s_1,u}-\mu_{s_2,u})\,du+\int_{s_2}^{s_1}(\mu_{u,u}-\mu_{s_2,u})\,du=
\int_0^t\int_{s_2}^{s_1}\Psi_t(x,y)(t-(x\vee y)) \, dx\, dy.
\label{drrddr}\end{equation} Comparing  \eqref{drrddr} with
\eqref{ige}, we see that condition \eqref{ige} is, at least
formally,  satisfied if $\Psi_t(x,y)$ has the form
\begin{equation} \Psi_t(x,y)= (t-(x\vee
y))^{-1}\int_{\R_+}(e^{-\tau\varkappa(x,y)(t-(x\vee
y))}-1+\tau\varkappa(x,y)(t-(x\vee
y)))\,\sigma(d\tau).\label{iuguzfv}\end{equation} To show that
this inserted into \eqref{ziguczug} gives indeed a solution of
\eqref{ige},  we have to verify that the $\Psi_t(x,y)$ given by
\eqref{iuguzfv} is integrable on $[0,t]^2$. Analogously to
\eqref{ibib}, we get \begin{align} &\int_0^t
\int_0^t\int_{(0,1]}|(t-(x\vee y))^{-1}
(e^{-\tau\varkappa(x,y)(t-(x\vee
y))}-1+\tau\varkappa(x,y)(t-(x\vee y)))| \,\sigma(d\tau)\, dx\,
dy\notag\\ &\qquad \le \int_0^t\int_0^t \int_{(0,1]}
\sum_{n=2}^\infty \frac{\tau^n\varkappa(x,y)^{n}(t-(x\vee
y))^{n-1}}{n!} \,\sigma(d\tau)\, dx\, dy\notag\\ &\qquad\le t^3
C_t^2 e^{tC_t}
\int_{(0,1]}\tau^2\,\sigma(d\tau)<\infty.\label{ztz}
\end{align}
Next,
 \begin{align} &\int_0^t
\int_0^t\int_{(1,+\infty)}|(t-(x\vee y))^{-1}
(e^{-\tau\varkappa(x,y)(t-(x\vee
y))}-1+\tau\varkappa(x,y)(t-(x\vee y)))| \,\sigma(d\tau)\, dx\,
dy\notag\\ &\qquad \le
\int_0^t\int_0^t\int_{(1,+\infty)}|(t-(x\vee y))^{-1}
(e^{-\tau\varkappa(x,y)(t-(x\vee y))}-1)| \,\sigma(d\tau)\, dx\,
dy\notag\\&\qquad\quad+t^2 C_t
\int_{(1,+\infty)}\tau\,\sigma(d\tau)\notag\\ &\qquad\le
2t^2C_t\int_{(1,+\infty)}\tau\,\sigma(d\tau),\label{gzgt}
\end{align}
where we used the estimate: $1-e^{-\alpha}\le \alpha$ for all
$\alpha\ge0$. Thus,  by \eqref{ziguczug} and
\eqref{iuguzfv}--\eqref{gzgt},  statement (a) holds for
\begin{equation} \mu_{s,t}=\mu_{0,t}+\int_0^t\int_0^s \int_{\R_+}
(t-(x\vee y))^{-1} (e^{-\tau\varkappa(x,y)(t-(x\vee
y))}-1+\tau\varkappa(x,y)(t-(x\vee y))) \,\sigma(d\tau)\,dx\,dy.
\label{huasgu}\end{equation}

Let us now suppose that (a), or equivalently \eqref{ige}, holds.
We set $s_2=0$ and $s_1=s$. Then, \eqref{ige} takes the following
form: \begin{align}& \int_0^t \int _{0}^{s} \int_{\R_+}
(e^{-\tau\varkappa(x,y)(t-(x\vee y))}-1+
\tau\varkappa(x,y)(t-(x\vee y))) \,\sigma(d\tau)\, dx\, dy\notag\\
&\qquad=\int_{s}^t
(\mu_{s,u}-\mu_{0,u})\,du+\int_{0}^{s}(\mu_{u,u}-\mu_{0,u})\,du.\label{igagdz}
\end{align}
Differentiating \eqref{igagdz} in $t$ yields for $s\le t$:
\begin{equation}\mu_{s,t}-\mu_{0,t} = \int_0^t \int_0^s
\int_{\R_+}\tau\varkappa(x,y)(1-e^{-\tau\varkappa(x,y)(t-(x\vee
y))})\,\sigma(d\tau)\, dx\, dy. \label{huteri}\end{equation} That
the integral on the right hand side of \eqref{huteri} is finite
may be shown analogously to \eqref{ztz}, \eqref{gzgt}.

Since for the $\mu_{s,t}$ given by \eqref{huasgu}  statement (a)
holds, this $\mu_{s,t}$  also satisfies \eqref{huteri}. Therefore,
\begin{align*}& \int_0^t\int_0^s \int_{\R_+} (t-(x\vee y))^{-1}
(e^{-\tau\varkappa(x,y)(t-(x\vee
y))}-1+\tau\varkappa(x,y)(t-(x\vee y)))
\,\sigma(d\tau)\,dx\,dy\\&\qquad =\int_0^t \int_0^s
\int_{\R_+}\tau\varkappa(x,y)(1-e^{-\tau\varkappa(x,y)(t-(x\vee
y))})\,\sigma(d\tau)\, dx\, dy,\end{align*} which implies the
equivalence of (a) and (b).

Let us now show the equivalence of (b) and (c). Analogously to
Kennedy (1994), we conclude  that (c) is equivalent to
\begin{equation}
\int_{\D'}\exp\bigg(-\int_s^t(F_{u,u}-F_{s,u})\bigg)\, d\nu=1
\label{sdhz}\end{equation} for all $s,t\ge0$, $s\le t$.
Analogously to the above, we see that \eqref{sdhz} is, in turn,
equivalent to
\begin{align}&\int_0^t \int_s^t \int_{\R_+}(e^{-\tau\varkappa(x,y)(t-(x\vee y))}-1+
\tau\varkappa(x,y)(t-(x\vee y)))\,\sigma(d\tau)\, dx\, dy\notag\\
&\qquad =\int_{s}^t
(\mu_{u,u}-\mu_{s,u})\,du\label{igeawz}\end{align} for all
$s,t\ge0$, $s\le t$. Setting in \eqref{ige} $s_2=s$ and $s_1=t$,
we see that \eqref{igeawz} is a special case of \eqref{ige}, so
that (b) implies (c). To show the inverse conclusion, we follow
Kennedy (1994). Differentiating \eqref{igeawz} in $t$ yields
\begin{equation}  \int_0^t \int_s^t\int_{\R_+} \tau\varkappa(x,y)(1-e^{-\tau\varkappa(x,y)(t-(x\vee
y))})\,\sigma(d\tau)\, dx\,dy=\mu_{t,t}-\mu_{s,t}\label{twzcf}
\end{equation} for all $s,t\ge0$, $s\le t$. Setting $s=0$ in the
latter equation gives
\begin{equation}  \int_0^t \int_0^t\int_{\R_+}
\tau\varkappa(x,y)(1-e^{-\tau\varkappa(x,y)(t-(x\vee
y))})\,\sigma(d\tau)\,
dx\,dy=\mu_{t,t}-\mu_{0,t}.\label{sdbhif}\end{equation}
Subtracting \eqref{twzcf} from \eqref{sdbhif} implies (b).\quad
$\square$

\begin{corollary} Suppose that the L\'evy measure $\sigma$
additionally satisfies \begin{equation}\label{zcudgzu}
\la\tau\ra_\sigma{:=}\int_{\R_+}\tau\,\sigma(d
\tau)<\infty.\end{equation} Suppose that statement \rom{(a)} of
Theorem~\rom{\ref{th1}} holds and suppose that the initial term
structure $\{\mu_{0,t},\, t\ge 0\}$ satisfies \begin{equation}
\mu_{0,t}\ge \int_0^t\int_0^t\varkappa(x,y)\,dx\, dy\cdot
\la\tau\ra_\sigma,\qquad t\ge0. \label{ucasdtz}\end{equation} Then
the  forward rate process $\{F_{s,t},\, 0\le s\le t<\infty\}$
  and the spot rate process $\{R_t,\, t\ge 0\}$  take on
non-negative values $\nu$-a\rom.s\rom. \label{cor1}\end{corollary}

\noindent{\it Proof}. By \eqref{fwzc} and Theorem~\ref{th1}, we
get $$
F_{s,t}=\mu_{0,t}-\int_0^t\int_0^s\int_{\R_+}\tau\varkappa(x,y)
e^{-\tau\varkappa(x,y)(t-(x\vee y
))}\,\sigma(d\tau)\,dx\,dy+\widetilde X_{s,t},\qquad t\ge0,\ 0\le
s\le t,$$ where $$\widetilde X_{s,t}{:=}X_{s,t}+
\int_0^t\int_0^s\varkappa(x,y)\,dx\,dy\cdot \la\tau\ra_\sigma.$$

Under  condition \eqref{zcudgzu}, the measure $\nu$ is
concentrated on the set of all signed measures of the form
$\sum_{n=1}^\infty \tau _n\delta_{(x_n,y_n)}(dx\, dy)-\la \tau
\ra_\sigma\, dx\,dy$, where  $\delta_a$ denotes the Dirac measure
with mass at $a$, $\tau _n\in\operatorname{supp}\sigma$, $n\in\N$,
and  $\{(x_n,y_n) \}_{n=1}^\infty$ is a locally finite set in
$\R^2$, see, e.g., Lytvynov (2003). Therefore, by \eqref{zdf},
$\widetilde X_{s,t} $ takes on non-negative values $\nu$-a.s.
Furthermore, it follows from \eqref{ucasdtz} that $$
\mu_{0,t}-\int_0^t\int_0^s\int_{\R_+}\tau\varkappa(x,y)
e^{-\tau\varkappa(x,y)(t-(x\vee y ))}\,\sigma(d\tau)\,dx\,dy\ge
0,\qquad t\ge 0,\ 0\le s\le t,$$ from where the statement follows.
\quad $\square$

Let us consider two examples of a measure $\nu$ satisfying the
assumptions of  Theorem~\ref{th1} and Corollary~\ref{cor1}.

{\bf Example 1.} ({\it Poisson sheet}) We take as $\nu$ the
centered Poisson measure $\pi_z$ with intensity parameter $z>0$,
see, e.g., Hida (1970). The L\'evy measure $\sigma$ has now the
form $z\delta_1$. Thus, the Fourier transform of $\pi_z$ is given
by $$
\int_{\D'}e^{i\la\omega,\varphi\ra}\,\pi_z(d\omega)=\exp\bigg[
\int_{\R^2}(e^{i
\varphi(x,y)}-1-i\varphi(x,y))\,z\,dx\,dy\bigg],\qquad
\varphi\in\D. $$ We set $\varkappa(x,y)\equiv1$. Then, $X_{s,t}$
given by \eqref{zdf} with the underlying probability measure
$\nu=\pi_z$ is, by definition,  a Poisson sheet, and for each
fixed $t>0$, $\{X_{s,t},\, 0\le s\le t\}$ is a centered Poisson
process with intensity parameter $tz$.
 Statement (b) of Theorem~\ref{th1} now
reads as follows: $$
\mu_{s,t}=\mu_{0,t}+z\big((2-s)e^{s-t}-2e^{-t}-s+st\big).$$
Condition \eqref{ucasdtz} now means $\mu_{0,t}\ge zt^2$, $t\ge0$.

{\bf Example 2.} ({\it Gamma sheet}) We take as $\nu$ the centered
gamma measure $\gamma_z$ with intensity parameter $z>0$, see,
e.g., Lytvynov (2003). The L\'evy measure $\sigma$ on $\R_+$ has
the form $$ \sigma(d\tau)=\frac{e^{-\tau}}{\tau}\, z\, d\tau.$$
The Fourier transform of $\gamma_z$ may be written as follows: $$
\int_{\D'}e^{i\la
\omega,\varphi\ra}\,\gamma_z(d\omega)=\exp\bigg(-\int_{\R^2}\big(\log(1-i\varphi(x,y))+\varphi(x,y)\big)\,
z\,dx\,dy \bigg),\qquad \varphi\in\D,\ |\varphi|<1.$$ We set
$\varkappa(x,y)\equiv1$. Then, $X_{s,t}$  given by \eqref{zdf}
with the underlying probability measure  $\nu=\gamma_z$ is, by
definition,  a gamma sheet, and for each $t>0$, $\{X_{s,t},\, 0\le
s\le t\}$ is a centered gamma process with intensity parameter
$tz$. Statement (b) of Theorem~\ref{th1} now reads as follows:
$$\mu_{s,t}=\mu_{0,t}+z \bigg(
st+2s+2(1+t)\log\bigg(\frac{1+t-s}{1+t}\bigg)-s\log(1+t-s)
\bigg).$$ Condition \eqref{ucasdtz} means $\mu_{0,t}\ge zt^2$,
$t\ge0$.\vspace{3mm}

It is possible to construct a model of forward interest rates
which unifies the approach of Kennedy (1994) to modeling the
forward interest rate with our approach. Indeed, consider
$F_{s,t}$ in the form \begin{equation}\label{uirwe}
F_{s,t}=\mu_{s,t}+X_{s,t}+Y_{s,t},\qquad 0\le s\le
t,\end{equation} where $\mu_{s,t}$ and $X_{s,t}$ are as in formula
\eqref{fwzc} (thus, as in our approach) and $Y_{s,t}$ is a
centered continuous Gaussian random field that is independent of
$X_{u,v}$, $0\le u\le v<\infty$, and has covariance $$
\operatorname{Cov}(Y_{s_1,t_1},Y_{s_2,t_2})=c(s_1\wedge
s_2,t_1,t_2),\qquad 0\le s_i\le t_i,\ i=1,2,$$ with a function $c$
satisfying $c(0,t_1,t_2)\equiv 0$ (as in Kennedy's approach).

The following theorem may be easily proved by combining the proof
of Theorem~1.1 in Kennedy (1994) and the proof of
Theorem~\ref{th1}.

\begin{theorem}
Theorem~\rom{\ref{th1}} remains valid for the  forward rates
$\{F_{s,t},\,0\le s\le t<\infty\}$  given by \eqref{uirwe} if we
set the deterministic term $\mu_{s,t}$ in statement \rom{(b)} to
be $$ \mu_{s,t}=\mu_{0,t} + \int_0^t \int_0^s
\int_{\R_+}\tau\varkappa(x,y)(1-e^{-\tau\varkappa(x,y)(t-(x\vee
y))})\,\sigma(d\tau)\, dx\, dy +\int_0^t c(s\wedge u,u,t)\,du$$
for all $s,t\ge0$, $s\le t$.
\end{theorem}

\noindent{\Large\bf References}\vspace{3mm}

\noindent Adler, R.J.:  The geometry of random fields. Chichester:
John Wiley \&\ Sons 1981\vspace{3mm}

\noindent Berezansky, Yu.M., Sheftel, Z.G.,  Us, G.F.: Functional
analysis,  Vol. 2. Basel: Birkh\"auser Verlag 1996\vspace{3mm}

\noindent Bj\"ork, T.,   Christensen, B.J.:  Interest rate
dynamics and consistent forward rate curves.  Math.\ Finance {\bf
9}, 323--348 (1999)\vspace{3mm}

\noindent Bj\"ork, T.,  Di Masi, G., Kabanov, Yu.,  Runggaldier,
W.: Towards a general theory of bond markets. Finance and
Stochast.\ {\bf  1}, 141--174 (1997a)\vspace{3mm}

\noindent Bj\"ork, T.,  Kabanov, Yu.,   Runggaldier, W.:  Bond
market structure in the presence of marked point processes. Math.\
Finance {\bf 7}, 211--239 (1997b) \vspace{3mm}

\noindent Gel'fand, I.M.,  Vilenkin, N.Ya.: Generalized functions,
Vol.~4. Applications of harmonic analysis. New York: Academic
Press 1964\vspace{3mm}

\noindent Heath, D.C., Jarrow, R.A.,  Morton, A.: Bond pricing and
the term structure of interest rates: A discrete time
approximation. J. Financial Quant.\ Anal.\ {\bf 25}, 419--440
(1990) \vspace{3mm}

\noindent   Heath, D.C.,  Jarrow, R.A.,  Morton, A.: Bond pricing
and the term structure of interest rates: A new methodology for
contingent claims valuation. Econometrica {\bf  60}, 77--105
(1992)\vspace{3mm}

\noindent Hida, T.: Stationary stochastic processes. Princeton:
Princeton University Press 1979 \vspace{3mm}

\noindent Ho, T., Lee, S.: Term structure movements and pricing
interest rate contingent claims. J. Finance {\bf  1}, 1011--1029
(1986) \vspace{3mm}

\noindent  Kennedy, D.P.:  The term structure of interest rates as
a Gaussian random field. Math.\ Finance {\bf 4}, 247--258 (1994)
\vspace{3mm}

\noindent Kennedy, D.P.: Characterizing Gaussian models of the
term structure of interest rates.  Math.\ Finance {\bf 7},
107--118 (1997)\vspace{3mm}

\noindent Kusuoka, S.: Term structure and SPDE. In:  S. Kusuoka
and T. Maruyama (eds.): Advances in Mathematical Economics,
Vol.~2. Tokyo: Springer 2000, pp.~67--85\vspace{3mm}

\noindent  L\"utkebohmert, E.: Endlich dimensionale Darstellungen
f\"ur das erweiterte Zinsmodell von Heath, Jarrow und Morton.
Diploma Thesis, Bonn: Bonn University 2002\vspace{3mm}

\noindent Lytvynov, E.:  Orthogonal decompositions for L\'evy
processes with an application to the gamma, Pascal, and Meixner
processes. Infin.\ Dimens.\ Anal.\ Quantum Probab.\ Relat.\ Top.\
{\bf 6}, 73--102 (2003)\vspace{3mm}

\noindent Mahnig, A.:   Modellierung der Zinsstrukturen durch ein
Poisson Sheet. Diploma Thesis, Bonn: Bonn University
2002\vspace{3mm}

\noindent  Popovici, S.A.:  Modellierung von Zinsstrukturkurven
mit Hilfe von stochastischen partiellen Differentialgleichungen.
Diploma Thesis, Bonn: Bonn University 2001\vspace{3mm}

\noindent Yalovenko, I.:  Modellierung des Finanzmarktes und
unendlich dimensionale stochastische Prozesse. Diploma Thesis,
Bochum: Bochum University 1998

\end{document}